\documentclass{article}   	
\usepackage{amsmath,amssymb,bm}   
\usepackage[hyphens]{url}
\usepackage[colorlinks,citecolor=blue,urlcolor=blue]{hyperref}

\newenvironment{newreferences}
               {\section*{References}
                \begin{list}{}{\setlength{\itemsep}{0pt}
                               \setlength{\parsep}{0pt}
                               \setlength{\labelwidth}{0pt}
                               \setlength{\leftmargin}{12pt}
                               \setlength{\labelsep}{0pt}}
                \setlength{\itemindent}{-12pt}
               }{\end{list}}

\title{Bertrand's analysis of baccarat}
\author{Stewart N. Ethier}
\date{}

\begin{document}

\maketitle

\begin{abstract}
Joseph Bertrand [1822--1900], who is often credited with a model of duopoly that has a unique Nash equilibrium, made another significant contribution to game theory.  Specifically, his 1888 analysis of baccarat was the starting point for Borel's investigation of strategic games in the 1920s.  In this paper we show, with near certainty, that Bertrand's results on baccarat were borrowed, without attribution, from an 1881 paper of Albert Badoureau [1853--1923].  In addition, we discuss Borel's criticisms of Bertrand's analysis, one of which helps to explain why Badoureau's contribution was overlooked until now.\medskip

\noindent Keywords: baccarat chemin de fer, noncooperative game, pure strategy, best response.
\end{abstract}

\section{Introduction}

In \textit{The History of Game Theory, Volume I: From the Beginnings to 1945} (Dimand and Dimand, 1996, p.~132), we find, ``As Borel (1924, p.~101) noted, the starting-point for his investigation of strategic games was the analysis of baccarat by Bertrand.''  In Borel's (1924) own words,\footnote{This is an edited version of Savage's (1953) translation.}
\begin{quote}
The only author who has studied a particular case of the problem we envisage is Joseph Bertrand, in the passage of his calculus of probabilities that he devotes to the draw at five in the game of baccarat; there he clearly distinguishes the purely mathematical side of the problem from the psychological side, because he asks, on the one hand, if the punter has an advantage in drawing at five when the banker knows the punter's manner of play and, on the other hand, if the punter has an advantage in drawing at five by letting the banker believe, if he can, that it is not his custom; he also poses the same questions for the non-draw at five. But, as we will see from what follows, this study is extremely incomplete; on the one hand, Bertrand does not investigate what would happen if the punter would draw at five in a certain fraction of the total number of coups and, on the other hand, he does not investigate whether the banker's opinion on the draw at five of the punter does not affect the way the banker plays in some cases where the punter does not have five, which may result in an advantage for either player in the whole game.  We see, by these brief indications, how complex is the problem approached by Bertrand, in spite of the rudimentary simplicity of the game of baccarat; [\dots].
\end{quote}

To describe Bertrand's results, we must first explain the rules of the two-person zero-sum game of baccarat (more precisely, baccarat chemin de fer). 

The game is dealt from a shoe comprising six 52-card decks mixed together. Denominations A, 2--9, 10, J, Q, K have values 1, 2--9, 0, 0, 0, 0, respectively, and suits are irrelevant. The total of a hand, comprising two or three cards, is the sum of the values of the cards, modulo 10. In other words, only the final digit of the sum is used to evaluate a hand. Two cards are dealt face down to Player and two face down to Banker, and each looks only at his own hand. The object of the game is to have the higher total (closer to 9) at the end of play. A two-card total of 8 or 9 is a \textit{natural}. If either hand is a natural, the game is over. If neither hand is a natural, Player then has the option of drawing a third card. If he exercises this option, his third card is dealt face up. Next, Banker, observing Player's third card, if any, has the option of drawing a third card. This completes the game, and the higher total wins. Winning bets on Player's hand are paid at even odds, with Banker, as the name suggests, playing the role of the bank. Losing bets on Player's hand are collected by Banker. Hands of equal total result in a tie (no money is exchanged). Since multiple players can bet on Player's hand, Player's strategy is restricted. He must draw on a two-card total of 4 or less and stand on a two-card total of 6 or 7. When his two-card total is 5, he is free to stand or draw as he chooses. (The decision is usually made by the player with the largest bet.)  In much of the French baccarat literature, a Player who always draws at 5 is called a \textit{tireur};  likewise, a Player who always stands at 5 is called a \textit{non-tireur}.  Banker, on whose hand no one can bet, has no constraints on his strategy in the classical version of the game considered here.

Bertrand (1888, pp.~39--40) states his problem precisely:

\begin{quote}
We must solve four problems:

The player having 5 and not requesting a card, what is the probability for him to win and what is that of a tie, when the banker, unaware that he has the point 5, knows that he has the custom, when he has this point, of not requesting a card?

The player having 5 and not requesting a card, what are the probabilities for him to win or to tie, when the banker, unaware that he has the point 5, knows that he has the custom, when he has this point, of requesting a card?

The player having 5 and requesting a card, what are the probabilities for him to win or to tie, when the banker knows his custom of requesting a card in this circumstance or when he thinks he knows that he does not request one?
\end{quote}

Keep in mind that, in the late 19th century, game theory did not yet exist, so terms such as ``pure strategy'' and ``best response,'' which we use below, were not yet in common usage. 

In the next section we show how to solve these four problems.  In Section~\ref{predecessors} we compare our results with those of Bertrand and those of Bertrand's predecessors in the analysis of baccarat, Dormoy (1872) and Badoureau (1881).  Finally, in Section~\ref{Borel} we discuss Borel's criticisms of Bertrand's analysis.

The principal characters in the story include:
\begin{itemize}
\item Joseph Bertrand [1822--1900], professor of mathematics at \'Ecole Polytechnique. He contributed to several areas of mathematics, including probability, in which he is best known for the ballot theorem and ``Bertrand's paradox'' (concerning a random chord of a circle).

\item \'Emile Dormoy [1829--1891], a French mining engineer and actuary.  He published numerous books, including three on card games (baccarat, bouillotte, and \'ecart\'e).  His best-known work is the 1878 book, \textit{Th\'eorie math\'e\-matique des assurances sur la vie}, in two volumes.

\item Albert Badoureau [1853--1923], a French mining engineer and amateur mathematician.  He is best known for his 1878 discovery of 37 of the 75 non-prismatic uniform polyhedra. He was also scientific advisor to Jules Verne and wrote a technical appendix to \textit{Topsy Turvy}.

\item \'Emile Borel [1871--1956], professor of mathematics at the Sorbonne and one of the great French mathematicians.  He was one of the founders of the theory of functions of a real variable, and he was an early contributor to game theory.
\end{itemize}

Despite their 49-year age difference, Borel knew Bertrand personally (Bru, 2006, note 12); in fact, Borel's wife was granddaughter of Bertrand's elder brother.  Moreover, it is very likely that Badoureau and Bertrand also knew each other.  Badoureau, having been ranked first in his class at \'Ecole Polytechnique, was responsible for writing the course handouts for his fellow students.  This would have included the first-year analysis course, which was offered in alternate years by Bertrand and Charles Hermite [1822--1901].  Records show that Bertrand taught the 1878 class and Badoureau was in the 1874 class.  This suggests that Badoureau's professor for this subject was Bertrand, not Hermite  (Bernard Bru, personal communication, 2021; Jacques Crovisier, personal communication, 2021).  It is also possible that Dormoy and Bertrand knew each other, for Dormoy ranked fourth in his class at \'Ecole Polytechnique, though this was while Bertrand was still a tutor in analysis, before being promoted to professor (Bru and Jongmans, 2001).

\section{Solution of Bertrand's problems}

There are four scenarios.  Player has a two-card total of 5 and can stand or draw.  Banker replies with a best response (unaware of Player's two-card total), assuming either that Player is a \textit{non-tireur} or that Player is a \textit{tireur}.  So the first step is to derive Banker's best response to each of Player's two pure strategies.  We assume, as is conventional, that cards are drawn with replacement (or, equivalently, that the shoe has infinitely many decks).
We can summarize these best responses in terms of two $8\times11$ incidence matrices,
\begin{equation*}
\arraycolsep=.7mm
\bm D_0=\left(
\begin{array}{ccccccccccc}
 1 & 1 & 1 & 1 & 1 & 1 & 1 & 1 & 1 & 1 & 1 \\
 1 & 1 & 1 & 1 & 1 & 1 & 1 & 1 & 1 & 1 & 1 \\
 1 & 1 & 1 & 1 & 1 & 1 & 1 & 1 & 1 & 1 & 1 \\
 1 & 1 & 1 & 1 & 1 & 1 & 1 & 1 & 0 & 0 & 1 \\
 0 & 0 & 1 & 1 & 1 & 1 & 1 & 1 & 0 & 0 & 1 \\
 0 & 0 & 0 & 0 & 0 & 1 & 1 & 1 & 0 & 0 & 1 \\
 0 & 0 & 0 & 0 & 0 & 0 & 1 & 1 & 0 & 0 & 0 \\
 0 & 0 & 0 & 0 & 0 & 0 & 0 & 0 & 0 & 0 & 0 \\
\end{array}
\right),\qquad
\bm D_1=\left(
\begin{array}{ccccccccccc}
 1 & 1 & 1 & 1 & 1 & 1 & 1 & 1 & 1 & 1 & 1 \\
 1 & 1 & 1 & 1 & 1 & 1 & 1 & 1 & 1 & 1 & 1 \\
 1 & 1 & 1 & 1 & 1 & 1 & 1 & 1 & 1 & 1 & 1 \\
 1 & 1 & 1 & 1 & 1 & 1 & 1 & 1 & 0 & 1 & 1 \\
 0 & 1 & 1 & 1 & 1 & 1 & 1 & 1 & 0 & 0 & 1 \\
 0 & 0 & 0 & 0 & 1 & 1 & 1 & 1 & 0 & 0 & 1 \\
 0 & 0 & 0 & 0 & 0 & 0 & 1 & 1 & 0 & 0 & 1 \\
 0 & 0 & 0 & 0 & 0 & 0 & 0 & 0 & 0 & 0 & 0 \\
\end{array}
\right),
\end{equation*}
in which rows are labeled by Banker's two-card total, 0--7, and columns are labeled by Player's third card, 0--9 and $\varnothing$, where $\varnothing$ means that Player stands.  Entries are 0 for stand and 1 for draw, with $\bm D_0$ giving Banker's best response when Player is known to be a \textit{non-tireur}, and $\bm D_1$ giving Banker's best response when Player is known to be a \textit{tireur}.

Notice that $\bm D_0$ and $\bm D_1$ differ at four entries:  At $(3,9)$, $(4,1)$, $(5,4)$, and $(6,\varnothing)$, $\bm D_0$ has a 0 and $\bm D_1$ has a 1.

How are $\bm D_0$ and $\bm D_1$ derived?  Let $p(i):=(1+3\delta_{i,0})/13$ ($i=0,1,\ldots,9$) be the distribution of the value of a third card, and $q(j):=(16+9\delta_{j,0})/(13)^2$ ($j=0,1,\ldots,9$) be the distribution of the value of a two-card hand.  Define the function $M$ on the set of nonnegative integers by $M(i):=\text{mod}(i,10)$.  Then, for $u\in\{0,1\}$, the $(j,k)$ entry of $\bm D_u$ is 1 if, when Banker's two-card total is $j$ and Player's third card is $k$, Banker's expected profit by drawing is greater than Banker's expected profit by standing, that is, if
$b_{u,1}(j,k)>b_{u,0}(j,k)$, where
\begin{equation}\label{b(j,k)}
\begin{split}
b_{u,0}(j,k)&:=\sum_{i=0}^{4+u}\text{sgn}(j-M(i+k))q(i)\bigg/\sum_{i=0}^{4+u}q(i),\\
b_{u,1}(j,k)&:=\sum_{i=0}^{4+u}\sum_{l=0}^9\text{sgn}(M(j+l)-M(i+k))q(i)p(l)\bigg/\sum_{i=0}^{4+u}q(i),
\end{split}
\end{equation}
when $j\in\{0,1,\ldots,7\}$ and $k\in\{0,1,\ldots,9\}$, and
\begin{equation}\label{b(j,e)}
\begin{split}
b_{u,0}(j,\varnothing)&:=\sum_{i=5+u}^7\text{sgn}(j-i)q(i)\bigg/\!\!\!\sum_{i=5+u}^7q(i),\\
b_{u,1}(j,\varnothing)&:=\sum_{i=5+u}^7\sum_{l=0}^9\text{sgn}(M(j+l)-i)q(i)p(l)\bigg/\!\!\!\sum_{i=5+u}^7q(i),
\end{split}
\end{equation}
when $j\in\{0,1,\ldots,7\}$.

One can go on to express the $2\times2^{88}$ payoff matrix for baccarat in terms of the conditional expectations \eqref{b(j,k)} and \eqref{b(j,e)} and then eliminate strictly dominated pure strategies and solve the game.  See Ethier (2010, Chapter 5) for full details.

In each of the equations in \eqref{b(j,k)} and \eqref{b(j,e)}, Player's two-card total is random, so these quantities will not help us solve Bertrand's problems, in which Player's two-card total is 5.  But we do need $\bm D_0$ and $\bm D_1$ to address these problems.  Let $f$ be a function on the set of integers and define
\begin{align}\label{A0v}
A_{0,v}&:=\sum_{\substack{0\le j\le7: \\ D_v(j,\varnothing)=1}}\sum_{l=0}^9 f(5-M(j+l))q(j)p(l)\bigg/\sum_{j=0}^7 q(j)\nonumber\\
&\qquad{}+\sum_{\substack{0\le j\le7: \\ D_v(j,\varnothing)=0}} f(5-j)q(j)\bigg/\sum_{j=0}^7 q(j)
\end{align}
and
\begin{align}\label{A1v}
A_{1,v}&:=\sum_{\substack{0\le j\le7,\,0\le k\le9: \\ D_v(j,k)=1}}\;\sum_{l=0}^9 f(M(5+k)-M(j+l))q(j)p(k)p(l)\bigg/\sum_{j=0}^7 q(j)\nonumber\\
&\qquad{}+\sum_{\substack{0\le j\le7,\,0\le k\le9: \\ D_v(j,k)=0}} f(M(5+k)-j)q(j)p(k)\bigg/\sum_{j=0}^7 q(j).
\end{align}
Then, with $f$ equal to the indicator of the set of positive integers, $A_{u,v}$ is $W_{u,v}$, the probability of a Player win; with $f$ equal to the indicator of $\{0\}$, $A_{u,v}$ is $T_{u,v}$, the probability of a Player tie; and with $f(x):=\text{sgn}(x)$, $A_{u,v}$ is $E_{u,v}$, Player's expected profit.  In all cases, $u=0$ (resp., $u=1$) if Player has a two-card total of 5 (with Banker unaware of this) and stands (resp., draws), and $v=0$ (resp., $v=1$) if Banker assumes that Player is a \textit{non-tireur} (resp., \textit{tireur}) and makes a best response.  

Notice another distinction between \eqref{b(j,k)}--\eqref{b(j,e)} and \eqref{A0v}--\eqref{A1v}.  In the former, the $u$ variable signifies Player's pure strategy, and the $v$ variable signifies Banker's stand-or-draw decision based on the available information.  In the latter, the $u$ variable signifies Player's stand-or-draw decision with a two-card total of 5, and the $v$ variable signifies what Banker assumes is Player's pure strategy, to which Banker makes a best response.

Computations show that
\begin{align}\label{E00}\tag{5a}
W_{0,0}&=\frac{792}{1781},&\quad T_{0,0}&=\frac{153}{1781},&\quad E_{0,0}&=-\frac{44}{1781},\\ \label{E01}\tag{5b}
W_{0,1}&=\frac{872}{1781},&\quad T_{0,1}&=\frac{169}{1781},&\quad E_{0,1}&=\frac{132}{1781},\\ \label{E10}\tag{5c}
W_{1,0}&=\frac{10352}{23153},&\quad T_{1,0}&=\frac{2928}{23153},&\quad E_{1,0}&=\frac{479}{23153},\\ \label{E11}\tag{5d}
W_{1,1}&=\frac{10176}{23153},&\quad T_{1,1}&=\frac{2976}{23153},&\quad E_{1,1}&=\frac{175}{23153}.
\end{align}
Alternatively,
\begin{align}\label{E00dec}\tag{6a}
W_{0,0}&\approx0.444694,&\quad T_{0,0}&\approx0.085907,&\quad E_{0,0}&\approx-0.024705,\\ \label{E01dec}\tag{6b}
W_{0,1}&\approx0.489613,&\quad T_{0,1}&\approx0.094891,&\quad E_{0,1}&\approx0.074116,\\ \label{E10dec}\tag{6c}
W_{1,0}&\approx0.447113,&\quad T_{1,0}&\approx0.126463,&\quad E_{1,0}&\approx0.020688,\\ \label{E11dec}\tag{6d}
W_{1,1}&\approx0.439511,&\quad T_{1,1}&\approx0.128536,&\quad E_{1,1}&\approx0.007558.
\end{align}

Although the formulas \eqref{A0v} and \eqref{A1v} are not difficult to evaluate, they depend crucially on $\bm D_0$ and $\bm D_1$, whose evaluation was a daunting task prior to the computer age.\footnote{None of the 19th-century baccarat authors got $\bm D_0$ and $\bm D_1$ exactly right.  This includes Billard (1883), Bertez\`ene (1889, 1896, 1897), Hoffmann (1891a,b), and Laun (1891), in addition to Dormoy (1872), Badoureau (1881), and Bertrand (1888).  We can add 20th-century baccarat authors Savigny (1906), Lafaye and Krauss (1927), and Huitte (1928) to this list.  (To be fair, Lafaye and Krauss assumed a different model of card distribution.)  Thus, the first to get $\bm D_0$ and $\bm D_1$ exactly right were Lafrogne (1927), Le Myre (1935), and Boll (1936).  Jules Louis Henry Lafrogne [1867--1933] was a Rear Admiral in the French Navy.  As far as we know, his baccarat book is his only publication.  ``Georges Le Myre'' was said (in a 1937 review of his book) to be the pseudonym of a professor of mathematics well known in Paris, but his identity is unknown today.  Marcel Boll [1886--1971] was professor of Chemistry and Electricity at l'\'Ecole des Hautes \'Etudes Commerciales in Paris.}  This may explain why a solution of the game of baccarat, promised by von Neumann (1928), would have to wait for Kemeny and Snell (1957), but even they did not use a computer but rather ``one of the old calculators.''\footnote{From a 3 August 2000 email from J. Laurie Snell to R\'egis Deloche.} Kemeny, incidentally, would later co-develop the BASIC programming language.

\section{Bertrand and his baccarat predecessors}\label{predecessors}

Bertrand had at least two predecessors in baccarat analysis.  

Dormoy's (1872) article (and nearly identical book (1873)) was the first mathematical study of baccarat chemin de fer, though he referred to it as ``baccarat tournant'' (rotating baccarat).  He obtained the equivalent of $\bm D_0$ and $\bm D_1$, in the form of two tables, with only one error in each.  He argued that it is indifferent whether Banker stands or draws at $(5,4)$ when Player is known to be a \textit{non-tireur}.  This is not really an error, just a consequence of rounding expectations to two decimal places (the standing expectation is $-299/1157$, and the drawing expectation is $-300/1157$; both round to $-0.26$).  He also had an error in $\bm D_1$ at $(6,6)$.

Dormoy addressed three of Bertrand's four problems, at least in terms of expectations, but rounded them to two or three decimal places.  He obtained $E_{0,0}\approx-0.024$, $E_{1,0}\approx0.012$, and $E_{1,1}\approx0.02$ (Section 79); these numbers should have been $-0.025$, $0.021$, and $0.01$.  He provided enough detail to allow one to check that his methodology was sound, so the errors were only computational.   

The second mathematical study of baccarat chemin de fer was by Badoureau (1881), who was aware of Dormoy's work.  He had $\bm D_0$ exactly correct but had three errors in $\bm D_1$, at $(4,1)$, $(4,9)$, and $(6,6)$.  Badoureau expressed $\bm D_0$, for example, in terms of a list of rules for Banker:

\begin{quote}
\begin{enumerate}
\item Always draw at baccarat [i.e., zero], one or two;

\item Always stand at seven;

\item Before a player who stands, draw at 3, at 4 or at 5 and stand at 6;

\item Before a player who has taken a card, always draw at 3 unless one has given 8 or 9; always draw at 4, unless one has given 8, 9, 0 or 1; only draw at 5 if one has given 5, 6 or 7; only draw at 6 if one has given 6 or 7. 
\end{enumerate}
\end{quote}

Badoureau addressed all four of Bertrand's problems and expressed his answers as fractions:
\begin{align}\tag{7a}\label{Ba00}
W_{0,0}&=\frac{792}{1781},&\quad T_{0,0}&=\frac{153}{1781},&\quad E_{0,0}&=-\frac{44}{1781},\\ \tag{7b}\label{Ba01}
W_{0,1}&=\frac{872}{1781},&\quad T_{0,1}&=\frac{169}{1781},&\quad E_{0,1}&=\frac{132}{1781},\\ \tag{7c}\label{Ba10}
W_{1,0}&=\frac{10352}{23153},&\quad T_{1,0}&=\frac{2928}{23153},&\quad E_{1,0}&=\frac{479}{23153},\\ \tag{7d}\label{Ba11}
W_{1,1}&=\frac{10288}{23153},&\quad T_{1,1}&=\frac{2800}{23153},&\quad E_{1,1}&=\frac{223}{23153}.
\end{align}
Actually, instead of finding $E_{u,v}$, Badoureau evaluated what he called Player's \textit{chances}, $C_{u,v}:=W_{u,v}+\frac12 T_{u,v}$.  Then $E_{u,v}=2C_{u,v}-1$, so $E_{u,v}$ and $C_{u,v}$ are effectively equivalent.

Finally, Bertrand (1888), who did not reveal $\bm D_0$ or $\bm D_1$ and cited no prior work, expressed his results to six decimal places:
\begin{align}\tag{8a}\label{Be00}
W_{0,0}&\approx0.444694,&\quad T_{0,0}&\approx0.085907,&\quad E_{0,0}&\approx-0.024706,\\ \tag{8b}\label{Be01}
W_{0,1}&\approx0.489612,&\quad T_{0,1}&\approx0.094890,&\quad E_{0,1}&\approx0.074115,\\ \tag{8c}\label{Be10}
W_{1,0}&\approx0.447113,&\quad T_{1,0}&\approx0.126463,&\quad E_{1,0}&\approx0.020689,\\ \tag{8d}\label{Be11}
W_{1,1}&\approx0.444348,&\quad T_{1,1}&\approx0.120935,&\quad E_{1,1}&\approx0.009631.
\end{align}
Actually, Bertrand found $L_{u,v}=1-W_{u,v}-T_{u,v}$, the probability of a Player loss, instead of $E_{u,v}=W_{u,v}-L_{u,v}=2W_{u,v}+T_{u,v}-1$. 

We find that \eqref{Be00}--\eqref{Be10}, \eqref{Ba00}--\eqref{Ba10}, and \eqref{E00}--\eqref{E10} (equiv., \eqref{E00dec}--\eqref{E10dec}) match exactly, ignoring rounding errors in the sixth decimal place.  But \eqref{Be11} matches \eqref{Ba11}, which differs from \eqref{E11} (equiv., \eqref{E11dec}).  How did Badoureau arrive at $W_{1,1}=10288/23153$, for example?  When we recompute \eqref{A0v}--\eqref{A1v} for the appropriate $f$ with $\bm D_1$ changed in three places (at $(4,1)$, $(4,9)$, and $(6,6)$) to account for Badoureau's errors, we obtain Badoureau's numbers exactly.  

For Bertrand to have derived his results independently of Badoureau, he would have had to have made no errors in the 88 entries of $\bm D_0$ and exactly the same three errors as Badoureau made in the 88 entries of $\bm D_1$.  This seems extremely unlikely, and perhaps there is no need to quantify it.  We conclude that Bertrand borrowed Badoureau's results without confirming them and without attribution.  As Sheynin (1994) pointed out, Bertrand had the custom of not fully acknowledging prior work.  

Was Bertrand's interpretation of the results more insightful than Badou\-reau's?  Bertrand (1888, p.~42) wrote:
\begin{quote}
If, without trying to outwit the other, from the beginning of the game he frankly declares his customs, he should draw at 5.

If the conventions of the game allow deception, he should stand at 5, making the banker think, if he can, that he has the custom of drawing.
\end{quote}
Compare this with Badoureau (1881):
\begin{quote}
If one were required to make known in advance the rule of conduct that one proposes to follow, it would be better to draw than to stand, [\dots].

The player must therefore stand at 5 while letting the banker believe the opposite.
\end{quote}
The conclusions are similar.  See Guilbaud (1961) for discussion.

\section{Borel's criticism of Bertrand's analyis}\label{Borel}

As we saw, the first reason that Borel characterized Bertrand's study of baccarat as ``extremely incomplete'' was its lack of consideration of mixed strategies.  But both Dormoy (1872) and, to a lesser extent, Badoureau (1881) did consider a Player $(\frac12,\frac12)$ mixed strategy, although not entirely correctly.  A Banker best response to this mixed strategy was regarded, before game theory, as the appropriate response when Banker was unaware of Player's customary behavior.  In effect, Dormoy found Banker's best move at $(j,k)\in\{0,1,\ldots,7\}\times\{0,1,\ldots,9,\varnothing\}$ by determining which of
$$
\frac{b_{0,v}(j,k)+b_{1,v}(j,k)}{2},\qquad v\in\{0,1\},
$$
is larger.  A more accurate approach would have compared the weighted averages
$$
\frac{89\,b_{0,v}(j,k)+105\,b_{1,v}(j,k)}{89+105},\qquad v\in\{0,1\},
$$
if $k\ne\varnothing$, or
$$
\frac{48\,b_{0,v}(j,\varnothing)+32\,b_{1,v}(j,\varnothing)}{48+32},\qquad v\in\{0,1\},
$$
where the factor 89 (resp., 105, 48, and 32), when divided by 137, is the conditional probability that Player's two-card total belongs to $\{0,1,2,3,4\}$ (resp., $\{0,1,2,3,4,5\}$, $\{5,6,7\}$, and $\{6,7\}$), given that it is neither 8 nor 9.

Badoureau, on the other hand, made a more egregious error by using in effect the averages
$$
\frac{E_{0,0}+E_{0,1}}{2}=\frac{44}{1781}\approx0.024705\;\;\text{and}\;\; \frac{E_{1,0}+E_{1,1}}{2}=\frac{351}{23153}\approx0.015160
$$
(or $327/23153\approx0.014123$ for the second expression if the error in $E_{1,1}$ is corrected).  This mistake was said by Boll (1936, p.~234) to be common:
\begin{quote}
The first idea that would come to mind would be that the advantage would be the average of the advantages [\dots].

This simplistic idea, which is accepted without reflection by most players, is \textit{completely wrong}: it must take into account all the details of the tactic.
\end{quote}
Indeed, one should use \eqref{A0v}--\eqref{A1v} with $f(x):=\text{sgn}(x)$ and $\bm D_v$ replaced by $\bm D_{1/2}$, defined to be equal to $\bm D_0$ except for 1s at $(3,9)$ and $(5,4)$.  The latter approach yields
$$
E_{0,1/2}=-\frac{44}{1781}\approx-0.024705\quad\text{and}\quad E_{1,1/2}=\frac{287}{23153}\approx0.012396
$$
for Player's standing and drawing expectations in this scenario.

With these numbers, Badoureau might not have come to the conclusion that ``it is almost indifferent to draw or to stand at 5.''  Incidentally, $\bm D_{1/2}$ is Banker's mandatory strategy in today's nonstrategic form of baccarat, often called baccarat punto banco.  We have argued elsewhere (Ethier and Lee, 2015) that this is not a coincidence.

It appears that Lasker (1929) was the first author to consider a mixed strategy other than $(\frac12,\frac12)$ in the baccarat setting, and his analysis was correct, albeit for a simplified version of the game.  See Bewersdorff (2020) for discussion.

The second of Borel's criticisms of Bertrand's study is arguably even more fundamental.  Can one answer the question of whether Player should draw at 5 by examining only coups in which Player has a two-card total of 5, or must all Player two-card totals be included in the analysis?  We would argue for the latter approach.

For example, the first of Bertrand's four problems was the following:
\begin{quote}
The player having 5 and not requesting a card, what is the probability for him to win and what is that of a tie, when the banker, unaware that he has the point 5, knows that he has the custom, when he has this point, of not requesting a card?
\end{quote}

We believe that the problem would have been better formulated as follows:
\begin{quote}
What is the probability for the player to win and what is that of a tie, under the following two assumptions: (a) If the player has 5, he does not request a card.  (b) The banker knows that the player has the custom, when he has the point 5, of not requesting a card.
\end{quote}

As we have seen, the solution to the original problem is
\begin{equation}\tag{9}\label{Bertrand1}
W_{0,0}=\frac{792}{1781},\quad T_{0,0}=\frac{153}{1781},\quad E_{0,0}=-\frac{44}{1781},
\end{equation}
whereas the solution of the reformulated problem is
\begin{equation}\tag{10}\label{Bertrand1rev}
\overline{W}_{0,0}=\frac{2152648}{4826809},\quad \overline{T}_{0,0}=\frac{447337}{4826809},\quad \overline{E}_{0,0}=-\frac{74176}{4826809}.
\end{equation}
The distinction is clear:  Both sets of statistics assume that Player is a \textit{non-tireur} and  Banker knows it.  But $\overline{E}_{0,0}$ in \eqref{Bertrand1rev} is Player's mean profit from an arbitrary coup, whereas $E_{0,0}$ in \eqref{Bertrand1} is Player's mean profit on only those coups in which he starts with a two-card total of 5.  The latter statistic, having been taken out of context, is of dubious value.

Another way to look at this is to consider the $2\times2$ matrix game with payoff matrix
\begin{equation}\tag{11}\label{2x2,E}
\bordermatrix{& \text{B: assumes P stands at 5} & \text{B: assumes P draws at 5} \cr
\text{P: stands at 5} & E_{0,0}&E_{0,1}\cr
\text{P: draws at 5} & E_{1,0}&E_{1,1}}
\end{equation}
Player is the row player, has a two-card total of 5, and stands or draws.  Banker is the column player, is unaware that Player has a two-card total of 5, and makes a best response to one of Player's two pure strategies.  The equilibrium for this game is not meaningful because after a few plays Banker will realize that Player always has a two-card total of 5 and adjust his strategy accordingly.  

On the other hand, if \eqref{2x2,E} were replaced by
\begin{equation*}
\bordermatrix{& \text{B: assumes P stands at 5} & \text{B: assumes P draws at 5} \cr
\text{P: stands at 5} & \overline{E}_{0,0}&\overline{E}_{0,1}\cr
\text{P: draws at 5} & \overline{E}_{1,0}&\overline{E}_{1,1}},
\end{equation*}
with Player not assumed to have a two-card total of 5, the game could be played repeatedly and the equilibrium would be meaningful, even if not exactly the game-theoretic solution of baccarat.  It is, in fact, essentially the solution found by Kendall and Murchland (1964), who were unaware of Kemeny and Snell (1957).

It is for this reason, we suspect, that, with one exception,\footnote{Lafaye and Krauss (1927, Chap.~XI) correctly solved two of Bertrand's four problems, including the one he got wrong.  Yet they stated, not entirely accurately,  that their results were in conformity with those of Bertrand.} no subsequent researcher attempted to reproduce Bertrand's numbers, and therefore no one noticed the identical errors in Bertrand's and Badoureau's work.

Above all, we want to emphasize that Borel (1924) was mistaken in writing, ``The only author who has studied a particular case of the problem we envisage is Joseph Bertrand, [\dots].''  Game theorists might point to Charles Waldegrave, who solved the game of le her in 1713, as a counter-example (Guilbaud, 1961; Dimand and Dimand, 1996, Chapter 7; Bellhouse and Fillion, 2015), but Albert Badoureau is a more contemporaneous counter-example.  Badoureau died in 1923 so did not have the opportunity to correct the record.  Therefore, nearly a century later, we do so here.

\section*{Acknowledgments}

We thank Bernard Bru, Jacques Crovisier, R\'egis Deloche, Persi Diaconis, and Don Schlesinger for their help.

\begin{newreferences}

\item Badoureau, Albert (1881) \'Etude sur le jeu de baccarat. \textit{La Revue scientifique de la France et de l'\'etranger} \textbf{1} (19 f\'evrier) 239--246.

\item Bellhouse, David R. and Fillion, Nicolas (2015) Le her and other problems in probability discussed by Bernoulli, Montmort and Waldegrave. \textit{Statistical Science} \textbf{30} (1) 26--39.

\item Bertez\`ene, Alfred (1889) \textit{Baccarat, roulette, Bourse, le jeu, des s\'eries, de l'\'equil\-ibre, la mort lente}. E. Dentu, Paris.


\item Bertez\`ene, Alfred (1896) \textit{Le baccarat: c'est le vol.} Librairie de la Voix de Paris, Paris.

\item Bertez\`ene, Alfred (1897) \textit{M\'emoire \`a l'Acad\'emie des sciences}. Librairie de la Voix de Paris, Paris. 

\item Bertrand, Joseph (1888) \textit{Calcul des probabilit\'es}. Gauthier-Villars et Fils, Paris.

\item Bewersdorff, J\"org (2020) Emanuel Lasker and game theory.  In \textit{Emanuel Lasker, Volume II: Choices and Chances, Chess and Other Games of the Mind} (R. Forster, M. Negele, and R. Tischbierek, eds.). Exzelsior Verlag GmbH, Berlin, pp.~270--305.

\item Billard, Ludovic (1883) \textit{Br\'eviaire du baccara exp\'erimental}.  Chez l’auteur, Paris.

\item Boll, Marcel (1936) \textit{La chance et les jeux de hasard: Loterie, boule, roulettes, baccara, 30 \& 40, d\'es, bridge, poker, belote, \'ecart\'e, piquet, manille}.  Librairie Larousse, Paris.  (\'Edition revue, 1948.)

\item Borel, \'Emile (1924) Sur les jeux ou interviennent le hasard et l'habilet\'e des joueurs.  In \textit{\'El\'ements de la th\'eorie des probabilit\'es}. Librairie Scientifique, Hermann, Paris, pp.~204--224.  Reprinted in \textit{Revue d'\'economie politique} \textbf{69} (1959) (2) 139--167.  English translation by Savage, Leonard J.  On games that involve chance and the skill of the players.  \textit{Econometrica} \textbf{21} (1953) (1) 101--115.

\item Bru, Bernard (2006) Les le\c{c}ons de calcul des probabilit\'es de Joseph Bertrand. ``Les lois du hasard.''  \textit{Journal Electronique d'Histoire des Probabilit\'es et de la Statistique} \textbf{2} (2) 1--44.

\item Bru, Bernard and Jongmans, Fran\c cois (2001) Joseph Bertrand.  In \textit{Statisticians of the Centuries} (C. C. Heyde and E. Seneta, eds.). Springer-Verlag, New York, pp.~185--189.

\item Dimand, Mary Ann and Dimand, Robert W. (1996) \textit{The History of Game Theory, Volume I: From the Beginnings to 1945.} Routledge, Taylor \& Francis Group, London and New York.

\item Dormoy, \'Emile (1872) Th\'eorie math\'ematique des jeux de hasard. \textit{Journal des actuaires fran\c cais} \textbf{1} 120--146, 232--257.

\item Dormoy, \'Emile (1873) \textit{Th\'eorie math\'ematique du jeu de baccarat}. Armand Anger, Paris.

\item Ethier, Stewart N. (2010) \textit{The Doctrine of Chances: Probabilistic Aspects of Gambling}. Springer, Berlin and Heidelberg.

\item Ethier, Stewart N. and Lee, Jiyeon (2015) The evolution of the game of baccarat. \textit{Journal of Gambling Business and Economics} \textbf{9} (2) 1--13.

\item Guilbaud, G. Th.\ (1961) Faut-il jouer au plus fin? (Notes sur l'histoire de la théorie des jeux).  In \textit{La d\'ecision}, Colloques Internationaux du CNRS \textbf{515}. Centre national de la recherche scientifique \'editions, Paris, pp.~171--182.

\item Hoffmann, Professor Louis [pseudonym for Lewis, Angelo John] (1891a) \textit{The Cyclop{\ae}dia of Card and Table Games}. George Routledge and Sons, London.

\item Hoffmann, Professor Louis [pseudonym for Lewis, Angelo John] (1891b) \textit{Baccarat Fair and Foul, Being an Explanation of the Game, and a Warning against Its Dangers}. George Routledge and Sons, London.

\item Huitte, G. (1928) \textit{Le baccara: \`a la banque, chemin de fer,  r\`egles \& conseils}. Librairie Bernardin-B\'echet, Paris. 

\item Kemeny, John G. and Snell, J. Laurie (1957) Game-theoretic solution of baccarat. \textit{American Mathematical Monthly} \textbf{64} (7) 465--469.

\item Kendall, M. G. and Murchland, J. D. (1964) Statistical aspects of the legality of gambling.  \textit{Journal of the Royal Statistical Society, Series A} \textbf{127} (3) 359--391.

\item Lafaye, P. and Krauss, E. (1927) \textit{Nouvel expos\'e de la th\'eorie math\'ematique du jeu de baccara}. Papeteries Nouvelles Imp., Paris.

\item Lafrogne, Amiral Jules (1927) \textit{Calcul de l'avantage du banquier au jeu de baccara}. Gauthier-Villars et Cie, Paris.

\item Lasker, Emanuel (1929) \textit{Das verst\"andige Kartenspiel}. A. Scherl, Berlin.  

\item Laun [pseudonym for Delauney, Julien-F\'elix] (1891) \textit{Trait\'e th\'eorique et pratique du baccara, contenant la r\`egle du jeu, examen des probabilit\'es du tirage, r\`egle et th\'eorie du baccara en banque, grecs et tricheries}. Delarue, Paris.

\item Le Myre, Georges [pseudonym, real name unknown] (1935) \textit{Le baccara}. Hermann \& Cie, Paris.

\item Savigny, G. B. de [pseudonym for Renaudet, Georges-Benjamin] (1906) \textit{Les jeux de hasard.  Le Baccara, r\`egles compl\`etes des grands cercles.} Librairie des Publications Populaires, Paris. 

\item Sheynin, Oscar (1994) Bertrand's work on probability.  \textit{Archive for History of Exact Sciences} \textbf{48} (2) 155--199.

\item Von Neumann, John (1928)  Zur Theorie der Gesellschaftsspiele.  \textit{Mathematische Annalen} \textbf{100} (1) 295--320.

\end{newreferences}

\end{document}